\documentclass[final,leqno,onefignum]{siamltex704}

\usepackage{fancyhdr}

\usepackage{paralist}

\usepackage{ifthen}

\usepackage[usenames]{color}

\usepackage{epsfig,amsfonts,latexsym,amsmath,amssymb,hyperref}

\usepackage{tabularx}
\usepackage[T1]{fontenc}
\usepackage{times}
\input xy
\xyoption{all}

\newtheorem{remark}[theorem]{Remark}

\newcommand{\R}{{\mathbb R}}
\newcommand{\ZZ}{{\mathbb Z}}

\newcommand{\phig}{\varphi_\gamma}
\newcommand{\phitgnk}{\varphi_{m,\tilde\nu,k}}
\newcommand{\utgnk}{u_{m,\tilde\nu,k}}
\newcommand{\phignk}{\varphi_{m,\nu,k}}
\newcommand{\ugnk}{u_{m,\nu,k}}

\newcommand{\ii}{\mathrm{i}}
\newcommand{\dd}{\mathrm{d}}
\newcommand{\pdpd}[2]{\frac{\partial #1}{\partial #2}}
\newcommand{\pdpdo}[2]{\frac{d #1}{d #2}}

\def\ee{\mbox{e}}
\def\dd{\mbox{d}}

\newcommand{\centerPar}[1]{\xi'_{#1}}

\newcommand{\lengthPar}[1]{L'_{#1}}
\newcommand{\lengthPer}[1]{L''_{#1}}

\newcommand{\eoc}{$\Diamond$}

\title{Multi-scale discrete approximations of Fourier integral operators associated with canonical transformations and caustics~\thanks{The first three authors were partially supported by
    the National Science Foundation under grant CMG-1025259 and are
    grateful for the stimulating environment at the MSRI in Berkeley
    where this research was initiated in the Fall 2010.}}
\author{Maarten V. de Hoop\thanks{Department of Mathematics, Purdue
    University, West Lafayette, IN ({\tt
      mdehoop@math.purdue.edu})}\and Gunther Uhlmann\thanks{Department
    of Mathematics, University of Washington, Seattle, WA, and
    Department of Mathematics, UC Irvine, Irvine, CA ({\tt
      gunther@math.washington.edu})} \and Andr\'as
  Vasy\thanks{Department of Mathematics, Stanford University,
    Stanford, CA ({\tt andras@math.stanford.edu})} \and Herwig
  Wendt\thanks{Department of Mathematics, Purdue University, West
    Lafayette, IN; now at CNRS, IRIT UMR 5505, University of Toulouse, France
    ({\tt herwig.wendt@irit.fr})}}

\begin{document}

\maketitle
\thispagestyle{empty}

\begin{abstract}
We develop an algorithm for the computation of general Fourier
integral operators associated with canonical graphs. The algorithm is
based on dyadic parabolic decomposition using wave packets and enables
the discrete approximate evaluation of the action of such operators on
data in the presence of caustics. The procedure consists in the
construction of a universal operator representation through the
introduction of locally singularity-resolving diffeomorphisms,
enabling the application of wave packet driven computation, and in the
construction of the associated pseudo-differential joint-partition of
unity on the canonical graphs. We apply the method to a parametrix of
the wave equation in the vicinity of a cusp singularity.
\end{abstract}

\section{Introduction}

In this paper, we develop an algorithm for applying Fourier integral
operators associated with canonical graphs using wave packets. To
arrive at such an algorithm, we construct a universal oscillatory
integral representation of the kernels of these Fourier integral
operators by introducing singularity resolving diffeomorphisms where
caustics occur.  The universal representation is of the form such that
the algorithm based on the dyadic parabolic decomposition of phase
space previously developed by the authors applies
\cite{Andersson2011}. We refer to \cite{Bradie1993,Cand`es2003,Cand`es2007,Cand`es2009} for related
computational methods aiming at the evaluation of the action of
Fourier integral operators.

The algorithm comprises a geometrical component, bringing the local
representations in universal form, and a wave packet component which
yields the application of the local operators. Here, we develop the
geometrical component, which consists of the following steps. First we
determine the location of caustics on the canonical relation of the
Fourier integral operator.  For each point on a caustic we determine
the associated specific rank deficiency and construct an appropriate
diffeomorphism, resolving the caustic in open neighborhoods of this
point. We determine the (local) phase function of the composition of
the Fourier integral operator and the inverse of the diffeomorphism in
terms of universal coordinates and detect the largest set on which it
is defined. We evaluate the preimage of this set on the canonical
relation. We continue this procedure until the caustic is covered with
overlapping sets, associated with diffeomorphisms for the
corresponding rank deficiencies.  Then we repeat the steps for each
caustic and arrive at a collection of open sets covering the canonical
relation.

In the special case of Fourier integral operators corresponding to
parametrices of evolution equations, for isotropic media, an
alternative approach for obtaining solutions in the vicinity of
caustics, based on a re-decomposition strategy following a multi-product representation of
the propagator, has been proposed previously
\cite{Andersson2011,Kumano-go1978,Rousseau2006}.
Unlike multi-product representations, our construction does not involve a subdivision of the evolution parameter and yields a single-step computation. Moreover, it is valid
for the general class of Fourier integral operators associated with
canonical graphs, allowing for anisotropy. 

The complexity of the algorithm for general Fourier integral operators
as compared to the non-caustic case arises from switching, in the sets
covering a small neighborhood of the caustics, from a global to a
local algorithm using a pseudodifferential partition of unity.

As an application we present the computation of a parametrix of the
wave equation in a heterogeneous, isotropic setting for long-time
stepping in the presence of caustics.

\begin{figure}
\centering
\includegraphics[width=0.8\linewidth]{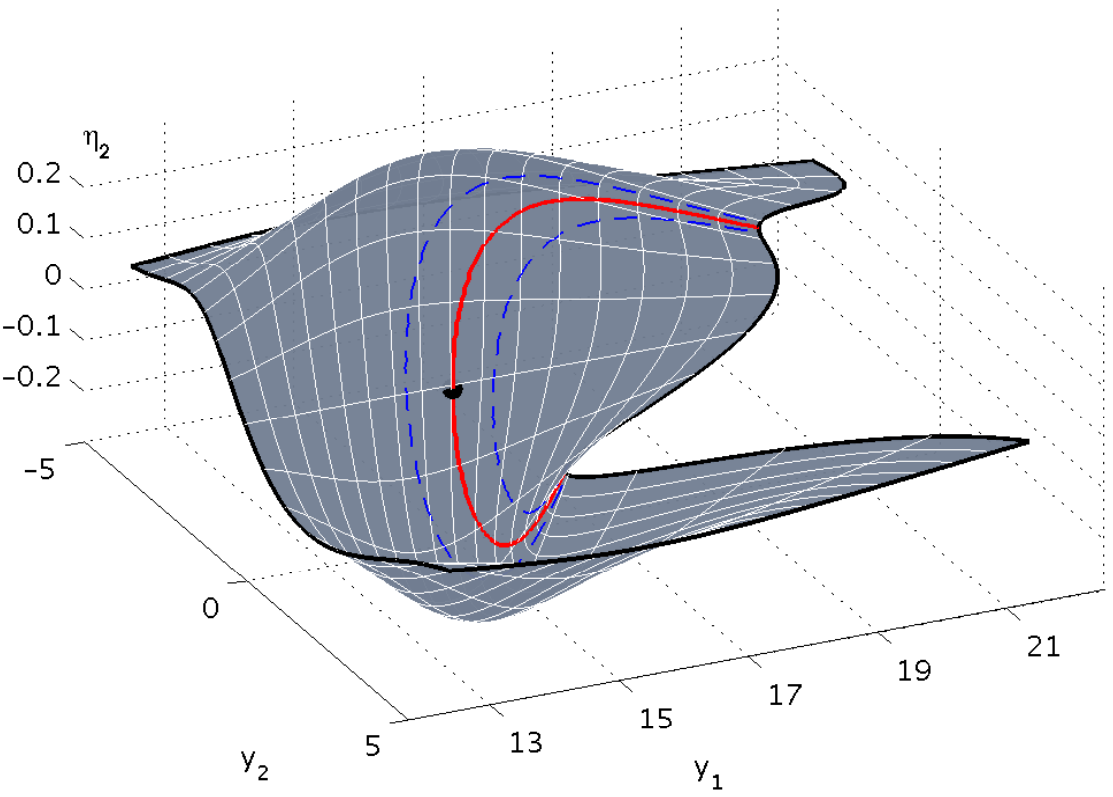}%
\caption{\label{fig:lagrangian} Projection $\bar\Lambda(y,\eta_2)$ of
  a slice $\xi=\xi_0$ of the canonical relation $\Lambda$ associated
  with a half-wave equation in the vicinity of a caustic (red solid
  line, the blue dashed lines indicate the neighborhood of the
  singularity) caused by a low velocity lens. The white solid lines
  are connected to a regular grid in $x$ by bi-characteristics. The
  black dot indicates the center of an open neighborhood of conjugate
  points $(x_0,\xi_0)\mapsto(y_0,\eta_0)$ for which the projection
  onto standard microlocal focal coordinates $(y,\xi)$ is not
  diffeomorphic.}
\end{figure}

\subsection*{Curvelets, wave packets}

We briefly discuss the (co)frame of curvelets and wave packets
\cite{Cand`es2006,Duchkov2009a,Smith1998}. Let $u \in L^2(\R^n)$ and
consider its Fourier transform, $\hat{u}(\xi) = \int u(x) \, \exp[-\ii
  \langle x,\xi \rangle ] \, \dd x$.

One begins with an overlapping covering of the positive $\xi_1$ axis
($\xi' = \xi_1$) by boxes of the form
\begin{equation}
\label{equ:box}
   B_k
   = \left[\centerPar{k} - \frac{\lengthPar{k}}{2},
           \centerPar{k} + \frac{\lengthPar{k}}{2}\right] \times
     \left[-\frac{\lengthPer{k}}{2},
            \frac{\lengthPer{k}}{2}\right]^{n-1} ,
\end{equation}
where the centers $\centerPar{k}$, as well as the side lengths
$\lengthPar{k}$ and $\lengthPer{k}$, satisfy the parabolic scaling
condition
\[
   \centerPar{k} \sim 2^k,\quad \lengthPar{k} \sim 2^k ,
   \quad \lengthPer{k} \sim 2^{k/2} ,\quad
                        \mbox{as $k \rightarrow \infty$} .
\]
Next, for each $k \ge 1$, let $\nu$ vary over a set of
$\sim 2^{k (n-1)/2}$ uniformly distributed unit vectors. 
Let $\Theta_{\nu,k}$ denote a choice of rotation matrix which maps $\nu$
to $e_1$, and
$
   B_{\nu,k} = \Theta_{\nu,k}^{-1} B_k .
$
In the (co-)frame construction, one
encounters two sequences of smooth functions on $\R^n$, $\hat\chi_{\nu,k}$ and
$\hat\beta_{\nu,k}$, each supported in $B_{\nu,k}$, so that
they form a co-partition of unity, 
$
   \hat{\chi}_0(\xi) \hat{\beta}_0(\xi) + \sum_{k \ge 1} \sum_{\nu}
   \hat{\chi}_{\nu,k}(\xi) \hat{\beta}_{\nu,k}(\xi) = 1 ,
$
and satisfy the estimates
$$
   |\langle \nu,\partial_{\xi} \rangle^j \, \partial_{\xi}^{\alpha}
    \hat{\chi}_{\nu,k}(\xi)| +
   |\langle \nu,\partial_{\xi} \rangle^j \, \partial_{\xi}^{\alpha}
    \hat{\beta}_{\nu,k}(\xi)|
   \le C_{j,\alpha} \, 2^{-k(j + |\alpha|/2)} .
$$
One then forms
$
   \hat{\psi}_{\nu,k}(\xi) = \rho_k^{-1/2}
             \hat{\beta}_{\nu,k}(\xi),\;
   \hat{\varphi}_{\nu,k}(\xi) = \rho_k^{-1/2}
              \hat{\chi}_{\nu,k}(\xi) ,
$
with $\rho_k=\textnormal{vol}(B_k)$, 
satisfying the
estimates
\begin{equation}
\label{eq:cdest}
\forall N:\quad
   \left.\begin{array}{l} |\varphi_{\nu,k}(x)| \\[0.3cm]
                |\psi_{\nu,k}(x)| \end{array} \right\}
          \le C_N 2^{k (n+1)/4} \, (\,
       2^k |\langle\nu,x\rangle| + 2^{k/2} \| x \| \,)^{-N}.
\end{equation}
To obtain a (co)frame, one introduces the integer
lattice: $X_m := (m_1,\dots,m_n)\in\ZZ^n$, the dilation matrix
$
   D_k = \frac{1}{2 \pi}
         \left(\begin{array}{lr} \lengthPar{k} & 0_{1 \times n-1} \\
     0_{n-1 \times 1} & \lengthPer{k} I_{n-1} \end{array}\right) ,
   \; \det \, D_k = (2\pi)^{-n} \rho_k ,
$
and points $x^{\nu,k}_m = \Theta_{\nu,k}^{-1} D_k^{-1} X_m$. The frame elements ($k \ge 1$) are then defined
in the Fourier domain as
$
   \hat{\varphi}_{\gamma}(\xi)
   = \hat{\varphi}_{\nu,k}(\xi) \,
         \exp[-\ii \langle x^{\nu,k}_m,\xi \rangle] ,\;
     \gamma = (m,\nu,k),
$
and similarly for $\hat{\psi}_{\gamma}(\xi)$. The function $\varphi_{\nu,k}$ is referred to as a wave packet. One obtains the transform pair
\begin{equation} \label{eq:ctrp}
   u_{\gamma} = \int u(x) \overline{\psi_{\gamma}(x)} \, \dd x ,
\quad\quad
   u(x) = \sum_{\gamma} u_{\gamma} \varphi_{\gamma}(x).
\end{equation}

\section{Fourier integral operators and caustics}

We consider Fourier integral operators, $F$, associated with canonical
graphs. We allow the formation of caustics.

\subsection{Oscillatory integrals, local coordinates}

Let $(y,x_{I_i},\xi_{J_i})$ be local coordinates on the canonical
relation, $\Lambda$ say, of $F$, and $S_i$ a corresponding generating
function: 
If, at a point on $\Lambda$, $(dy,dx_I)$ are linearly independent and $dx_J$ vanishes, then $(dy,dx_I,d\xi_J)$ are coordinates on $\Lambda$ nearby, $I\cup J=\{1,\dots,n\},\,I\cap J=\{\emptyset\}$, and one can parameterize $\Lambda$ as $\langle X_J(y,x_I,\xi_J)-x_J,\xi_J\rangle$, where $x_J=X_J(y,x_I,\xi_J)$ locally on $\Lambda$ (cf. \cite{Hormander1985}, Thm 21.2.18). The fact that a (possibly empty) set $I$ exists follows from the canonical graph property, i.e. that $(y,\eta)$ are local coordinates and $dy$ linearly independent.
Then
\begin{equation}
\begin{array}{rcrcl}
   \displaystyle{x_{J_i}
          = \frac{\partial S_i}{\partial \xi_{J_i}}} & , &
   \xi_{I_i} &=&
       \displaystyle{-\frac{\partial S_i}{\partial x_{I_i}}} ,
\\
   & & \eta &=& \displaystyle{\frac{\partial S_i}{\partial y}} .
\end{array}
\end{equation}
The coordinates are standardly defined on (overlapping) open sets
$O_i$ in $\Lambda$, that is, $(y,x_{I_i},\xi_{J_i}) \to
r(y,x_{I_i},\xi_{J_i})$ is defined as a diffeomorphism on $O_i$; let
$i = 1,\ldots,N$. The corresponding partition of unity is written as
\begin{equation}
   \sum_{i=1}^N \mathit{\Gamma}_i(r) = 1 ,\quad r \in \Lambda .
\end{equation}
In local coordinates, we introduce
\begin{equation}
   \bar{\mathit{\Gamma}}_i(y,x_{I_i},\xi_{J_i})
               = \mathit{\Gamma}_i(r(y,x_{I_i},\xi_{J_i})) .
\end{equation}
Then $(F \varphi_{\gamma})(y) = \sum_{i=1}^N
(F_i \varphi_{\gamma})(y)$ with
\begin{equation}
   (F_i \varphi_{\gamma})(y)= \int \! \int
           \bar{\mathit{\Gamma}}_i(y,x_{I_i},\xi_{J_i})
   a_i(y,x_{I_i},\xi_{J_i})
   \exp[\ii (\underbrace{S_i(y,x_{I_i},\xi_{J_i})
        - \langle \xi_{J_i},x_{J_i} \rangle}_{\phi(y,x,\xi_{J_i})}]\
   \varphi_{\gamma}(x)\ 
   \mathrm{d} x \mathrm{d}\xi_{J_i} .
\end{equation}
The amplitude $a_i(y,x_{I_i},\xi_{J_i})$ is complex and accounts for
the KMAH index. 

We let $\Sigma_{\phi}$ denote the stationary point set (in $\theta$)
of $\phi = \phi(y,x,\theta)$. The amplitude can be identified with a
half-density on $\Lambda$. One defines the $2n$-form $d_{\phi}$ on
$\Sigma_{\phi}$,
\[
   d_{\phi} \wedge d\left(\pdpd{\phi}{\theta_1}\right)
     \wedge \ldots \wedge d\left(\pdpd{\phi}{\theta_N}\right)
   = d y_1 \wedge \ldots \wedge d y_n \wedge
     d x_1 \wedge \ldots \wedge d x_n \wedge
     d\theta_1 \wedge \ldots \wedge d\theta_N .
\]
In the above, we choose $\lambda = (y,x_{I},\pdpd{\phi}{x_J})$ as
local coordinates on $\Lambda$, while $\theta = \xi_J$. Then we get
\[
   d_{\phi} = |\Delta_{\phi}|^{-1}
              |d\lambda_1 \wedge \ldots \wedge d\lambda_{2n}| ,\quad
   \Delta_{\phi} = \left| \begin{array}{cc}
            \frac{\partial^2 \phi}{\partial x_J \partial x_J} &
            \frac{\partial^2 \phi}{\partial \xi_J \partial x_J} \\
            \frac{\partial^2 \phi}{\partial x_J \partial \xi_J} &
            \frac{\partial^2 \phi}{\partial \xi_J \partial \xi_J}
                          \end{array} \right| = -1 ;
\]
$\lambda$ is identified with $(y,x_I,\xi_J)$. The corresponding
half-density equals $|\Delta_{\phi}|^{-1/2} |d\lambda_1 \wedge \ldots
\wedge d\lambda_{2n}|^{1/2}$.

Densities on a submanifold of the cotangent bundle are associated with
the determinant bundle of the cotangent bundle. Let $a_i^0$ denote the
leading order homogeneous part of $a_i$. The principal symbol of the
Fourier integral operator then defines a half-density, $a_i^0
d_{\phi}^{1/2}$. That is, for a change of local coordinates, if the
transformation rule for forms of maximal degree is the multiplication
by a Jacobian $\jmath$, then the transformation rule for a
half-density is the multiplication by $|\jmath|^{1/2}$. In our case,
of canonical graphs, we can dispose of the description in terms of
half-densities and restrict to zero-density amplitudes on $\Lambda$.

\subsection{Propagator}
\label{sec:prop}
The typical case of a Fourier integral operator associated with
a canonical graph is the parametrix for an evolution equation \cite{Duchkov2010,Duchkov2008},
\begin{equation}
\label{equ:evo:equ}
   [\partial_t + \ii P(t,x,D_x)] u(t,x) = 0 ,\quad u(t_0,x) = \phig(x)
\end{equation}
on a domain $X \subset \R^n$ and a time interval $[t_0,T]$, where
$P(t,x,D_x)$ is a pseudodifferential operator with symbol in
$S^1_{1,0}$; we let $p$ denote the principal symbol of $P$.

For every $(x,\xi) \in T^* X \backslash \{0\}$, the integral curves
$(y(x,\xi;t,t_0),\eta(x,\xi;t,t_0))$ of
\begin{equation}
\label{equ:ctsystem}
   \pdpdo{y}{t} = \pdpd{p(t,y,\eta)}{\eta} ,\quad
   \pdpdo{\eta}{t} = -\pdpd{p(t,y,\eta)}{y} ,
\end{equation}
with initial conditions $y(x,\xi;t_0,t_0) = x$ and
$\eta(x,\xi;t_0,t_0) = \xi$ define the transformation, $\chi$, from
$(x,\xi)$ to $(y,\eta)$, which generates the canonical relation of the
parameterix of \eqref{equ:evo:equ}, for a given time $t = T$; that is,
$(y(x,\xi),\eta(x,\xi)) = (y(x,\xi;T,t_0),\eta(x,\xi;T,t_0))$.

The perturbations of $(y,\eta)$ with respect to initial conditions
$(x,\xi)$ are collected in a propagator matrix,
\begin{equation}
\label{equ:fundmat}
   \Pi(x,\xi;t,t_0)
   = \left( \begin{array}{ c c }
     W_1 & W_2 \\
     W_3& W_4
   \end{array} \right)
   = \left( \begin{array}{ c c }
     \partial_x y & \partial_\xi y \\
     \partial_x \eta & \partial_\xi \eta
   \end{array} \right) ,
\end{equation} 
which is the solution to the $2n \times 2n$ system of differential
equations
\begin{equation}
\label{equ:perturbS}
\pdpdo{\Pi}{t}(x,\xi;t,t_0) = \left( \begin{array}{cc}
   \displaystyle{\frac{\partial^2 p}{
            \partial \eta \partial y}(t,y,\eta)} &
   \displaystyle{\frac{\partial^2 p}{
            \partial \eta \partial \eta}(t,y,\eta)} \\
   \displaystyle{-\frac{\partial^2 p}{
            \partial y \partial y}(t,y,\eta)} &
   \displaystyle{-\frac{\partial^2 p}{
            \partial y \partial \eta}(t,y,\eta)}
   \end{array} \right) \Pi(x,\xi;t,t_0) ,
\end{equation}
known as the Hamilton-Jacobi equations, supplemented with the initial
conditions \cite{Cerveny2001,Cerveny2007}
\begin{equation}
   \Pi(x,\xi;t_0,t_0) = \left( \begin{array}{cc} \mathbb{I}&{0}\\{0} & \mathbb{I}\end{array}\right).
\end{equation}
Away from caustics the generating function of $\Lambda$ is $S =
S(y,\xi)$ ($I_i = \emptyset$), which satisfies
\begin{eqnarray}
\label{equ:d2S1}
\pdpd{^2S}{y\partial\xi}(y,\xi)
     &=&\left.\pdpd{x}{y}\right|_\xi = W_1^{-1} , \\
\pdpd{^2S}{\xi^2}(y,\xi)
     &=&\left.\pdpd{x}{\xi}\right|_y=\left.-\pdpd{x}{y}\right|_\xi\left.\pdpd{y}{\xi}\right|_y = -W_1^{-1} W_2 , \\
\label{equ:d2S2}
\pdpd{^2S}{y^2}(y,\xi)
     &=&\left.\pdpd{\eta}{y}\right|_\xi=\left.\pdpd{\eta}{x}\right|_\xi\left.\pdpd{x}{y}\right|_\xi = W_3 W_1^{-1} ,
\end{eqnarray}
upon substituting $x = x(y,\xi;t_0,T))$ denoting the backward solution
to \eqref{equ:ctsystem} with initial time $T$, evaluated at $t_0$. The
leading-order amplitude follows to be
\begin{equation}
\label{equ:amp0}
   a(y,\xi / |\xi|) = \sqrt{1 /
      \det W_1(x(y,\xi / |\xi|;t_0,T),\xi / |\xi|;T,t_0)}  ,
\end{equation}
reflecting that $a$ is homogeneous of degree $0$ in $\xi$.

\medskip
In the vicinity of caustics, we need to choose different coordinates. Admissible coordinates are directly related to the possible rank
deficiency of $W_1$: One determines
the null space of the matrix $W_1$ and rotates the coordinates such that
the null space is spanned by the columns indexed by the set
$I_i$. Then $(y,x_{I_i},\xi_{J_i})$ form local coordinates on the
canonical relation $\Lambda$, as in the previous subsection, and $O_i$ is given by the set for which the columns indexed by  $I_i$ span the null space of $W_1$.

\section{Singularity resolving diffeomorphisms}

\begin{figure}
\centering
\includegraphics[width=0.8\linewidth]{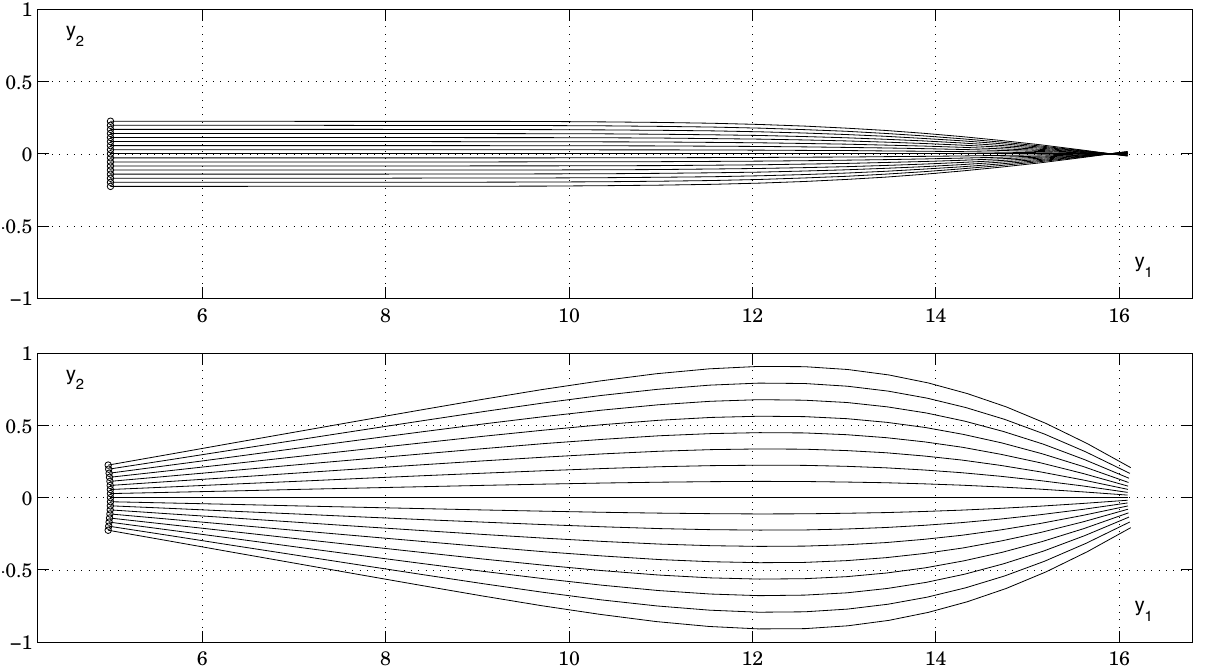}
\caption{\label{fig:canonical} Illustration of canonical relations $\chi$ (top) and $\check\chi$ (bottom) of operators $F$ and $\check F$ associated with a half-wave equation: (bi-)characteristics (''rays'') in $y$ for initial conditions $(x_2=x_{2,0},\;\xi=\xi_0)$ and  $(\tilde x_2=\tilde x_{2,0},\;\tilde\xi=\tilde\xi_0)$, respectively, for evolution through a low velocity lens (see Section \ref{sec:example}). The black circles on the left indicate the conjugate points corresponding to the initial conditions.}
\end{figure}

We consider the matrix $W_1(x(y,\xi;t_0,T),\xi;T,t_0)$ for given
$(T,t_0)$ at $y_0 = y(x_0,\xi_0;T,t_0)$ and $\xi = \xi_0$ and
determine its rank. Suppose it does not have full rank at this
point. We construct a diffeomorphism which removes this rank
deficiency in a neighborhood of $r_0 = (y_0,\eta_0;x_0,\xi_0) \in
\Lambda$, where $\eta_0 = \eta(x_0,\xi_0;T,t_0)$.

To be specific, we rotate coordinates, such that $\xi_0=(1,0,\dots,0)$
(upon normalization). Let us assume that the row associated with the
coordinate $x_2$ generates the rank
deficiency. 
(There could be more than one row / coordinate.)  We then
introduce the diffeomorphism,
\[
   Q : x \mapsto \tilde x
      = (x_1-\frac{\alpha}{2}(x_2-(x_0)_2)^2,x_2,\ldots,x_n) ;
\]
to preserve the symplectic form, we map
\[
   \xi \mapsto \tilde\xi
     = (\xi_1,\xi_2+\alpha(x_2-(x_0)_2) \, \xi_1,\xi_3,\ldots,\xi_n) ,
\]
yielding a canonical transformation $C_Q :\ (x,\xi) \mapsto (\tilde
x,\tilde\xi)$. We note that $C_Q(x_0,\xi_0) = (x_0,\xi_0)$. 

The canonical transformation, $C_Q^{-1}$, associated with $Q^{-1}$ is
given by
\begin{eqnarray*}
   \tilde x \to x &=& (\tilde x_1 + \frac{\alpha}{2}
         (\tilde x_2 - x_{0,2})^2,\tilde x_2,\ldots,\tilde x_n) ,
\\
   \tilde \xi \to \xi &=& (\tilde \xi_1,\tilde \xi_2 - \alpha
         (\tilde x_2 - x_{0,2}) \, \tilde\xi_1,\ldots,\tilde \xi_n) .
\end{eqnarray*}
We introduce the pull back, $Q^* u(\tilde x) = u(Q^{-1}(\tilde x)) =
u(\tilde x_1 +\frac{\alpha}{2}(\tilde x_2-(x_0)_2)^2,\tilde
x_2,\dots,\tilde x_n)$. 

\subsection{Fourier integral representations of $Q$ and $Q^{-1}$}
The diffeomorphism $Q$ can be written in the form of an invertible Fourier integral operator with unit amplitude and canonical relation given as the graph of $C_Q$.
To see this, we write $(Q^* u)(\tilde x) = u(X(\tilde x))$, $((Q^{-1})^* \tilde
u)(x) = \tilde u(\tilde X(x))$. That is, $X = Q^{-1}$ and $\tilde{X} =
Q$. The diffeomorphisms $Q$ and $Q^{-1}$ define the Fourier integral
operators with oscillatory integral kernels,
\begin{equation}
   A_Q(\tilde x,x) = 
   \int \ee^{-\ii \langle\xi,x-X(\tilde x)\rangle} \dd\xi ,
\quad
   A_{Q^{-1}}(x,\tilde x) =
   \int\ee^{-\ii \langle\tilde \xi,\tilde x-\tilde X( x)\rangle}
                  \dd\tilde \xi .
\end{equation}
The generating functions are
\[
   S_Q(\tilde x,\xi) = \langle\xi,X(\tilde x)\rangle ,
\quad
   S_{Q^{-1}}(x,\tilde \xi) = \langle\tilde \xi,\tilde X( x)\rangle ,
\]
respectively. The canonical relations are the graphs of $C_Q$ and
$C_{Q^{-1}}$, and are given by
\[
   \Lambda_{Q} = \{ (\tilde x = X^{-1}(x),
   \langle\xi,\partial_{\tilde x}X\rangle|_{
                     \tilde x = X^{-1}(x)};x,\xi) \} ,
\quad 
   \Lambda_{Q^{-1}} = \{ (x = \tilde X^{-1}(\tilde x),
   \langle\tilde \xi,\partial_{ x}\tilde X\rangle|_{
                x = \tilde X^{-1}(\tilde x)};\tilde x,\tilde \xi) \} .
\]
The Hessians yield a unit amplitude:
\[
   \left| \det \frac{\partial^2 \langle\xi,X(\tilde x)\rangle}{
               \partial \tilde x \partial \xi} \right| = 1 ,
\quad
   \left| \det \frac{\partial^2 \langle\tilde \xi,\tilde X( x)\rangle}{
               \partial x \partial \tilde \xi} \right| = 1 .
\]

Substituting the particular diffeomorphism, we obtain:
\begin{eqnarray*}
\partial_{ x}\tilde X|_{x = \tilde X^{-1}(\tilde x)}
&=&\left(\begin{array}{cccc}1&-\alpha(\tilde x_2-x_{0,2})&0&\cdots\\0&1&0&\cdots\\0&0&1&\cdots\\\vdots&\vdots&\vdots&\ddots\end{array}\right)\\
\langle\tilde \xi,\partial_{ x}\tilde X\rangle|_{x = \tilde X^{-1}(\tilde x)}
&=&\left(\begin{array}{c}\tilde\xi_1\\\tilde\xi_2-\alpha(\tilde x_2-x_{0,2})\tilde\xi_1\\\vdots\end{array}\right).
\end{eqnarray*}

The corresponding propagator matrices are hence given by
given by
\begin{eqnarray}
\Pi_Q &=&
\left(
\begin{array}{cc}
\pdpd{\tilde x}{x} & \pdpd{\tilde x}{\xi} \\
\pdpd{\tilde\xi}{x} & \pdpd{\tilde\xi}{\xi} 
\end{array}
\right)
=
\left(
\begin{array}{cccccccc}
1 & -\alpha(x_2 - x_{0,2}) & 0 & \cdots & 0 & 0 & 0  &\cdots \\
     0 & 1 & 0 & \cdots & 0 & 0 & 0  &\cdots \\
     0&0& 1 & \cdots & 0 & 0 & 0  &\cdots \\
     \vdots&\vdots&\vdots&\ddots & \vdots & \vdots & \vdots  &\ddots  \\
   0 & 0& 0 & \cdots & 1 & 0& 0 & \cdots \\
   0 & \alpha\xi_1& 0 & \cdots  & \alpha(x_2-x_{0,2}) & 1 & 0 & \cdots \\
     0&0& 0 & \cdots &0&0&1& \cdots \\
     \vdots&\vdots&\vdots&\ddots&\vdots&\vdots&\vdots&\ddots
\end{array}
\right) ,
\\
\Pi_Q^{-1}&=&
\left(
\begin{array}{cc}
\pdpd{x}{\tilde x} & \pdpd{x}{\tilde \xi} \\
\label{equ:propQ1}
\pdpd{\xi}{\tilde x} & \pdpd{\xi}{\tilde \xi} 
\end{array}
\right)
=
\left(
\begin{array}{cccccccc}
1 & \alpha(\tilde x_2 - x_{0,2}) & 0 & \cdots & 0 & 0 & 0  &\cdots \\
     0 & 1 & 0 & \cdots & 0 & 0 & 0  &\cdots \\
     0&0& 1 & \cdots & 0 & 0 & 0  &\cdots \\
     \vdots&\vdots&\vdots&\ddots & \vdots & \vdots & \vdots  &\ddots  \\
   0 & 0& 0 & \cdots & 1 & 0& 0 & \cdots \\
   0 & -\alpha\tilde\xi_1& 0 & \cdots  & - \alpha(\tilde x_2-x_{0,2}) & 1 & 0 & \cdots \\
     0&0& 0 & \cdots &0&0&1& \cdots \\
     \vdots&\vdots&\vdots&\ddots&\vdots&\vdots&\vdots&\ddots
\end{array}
\right) ,
\end{eqnarray}
which are easily verified to be symplectic matrices.  In the more
general case, each coordinate $x_j$ generating a rank deficiency
yields additional non-zero entry pairs $\pdpd{\tilde x_1}{x_j}$,
$\pdpd{x_1}{\tilde x_j}$, $\pdpd{\tilde \xi_j}{x_j}$,
$\pdpd{\xi_j}{\tilde x_j}$, and $\pdpd{\tilde \xi_j}{\xi_1}$,
$\pdpd{\xi_j}{\tilde\xi_1}$ in the above propagator matrices.

\subsection{Operator composition}
It follows that the composition $(\tilde x,
\tilde\xi)\stackrel{C_Q^{-1}}{\mapsto}(x,\xi)\stackrel{\chi
}{\mapsto}(y,\eta)$ generates the graph of a canonical transformation,
$\check{\chi}$ say, which can be parametrized by $(y,\tilde \xi)$
locally on an open neighborhood of $(y_0,\tilde{\xi}(x_0,\xi_0))$. We
denote the corresponding generating function by $\check S = \check
S(y,\tilde \xi)$. We can compose $F$ with $Q^{-1}$ as Fourier integral
operators: $\check F = F Q^{-1}$. The canonical relation of $\check F$
is the graph of $\check \chi$. In summary:

$$
\xymatrix@R=2.5pc@C=2.8pc{
(x,\xi) \ar@/_/[dd]_{Q:\,C_Q}\ar@{->}[rr]^{F:\,\chi}& & (y,\eta) \\
&\qquad\qquad&\\
(\tilde x,\tilde\xi)\ar@{->}[uurr]_{\check F = F Q^{-1}:\,\check\chi} 
\ar@/_/[uu]_{Q^{-1}:\,C_{Q^{-1}}}& & 
}
$$
For each given type of rank deficiency (here, in $x_2$) and each $(x_0,\xi_0)$ within this class, there is an open set $O_{(x_0,\xi_0)}$ on which the coordinates $(I,J)$ are valid. These sets form an open cover, and we obtain a family of diffeomorphisms parametrized by $(x_0,\xi_0)$; there exists a locally finite subcover, and we only need a discrete set to resolve the rank deficiencies everywhere.
%
\begin{figure}
\centering
\includegraphics[width=0.7\linewidth]{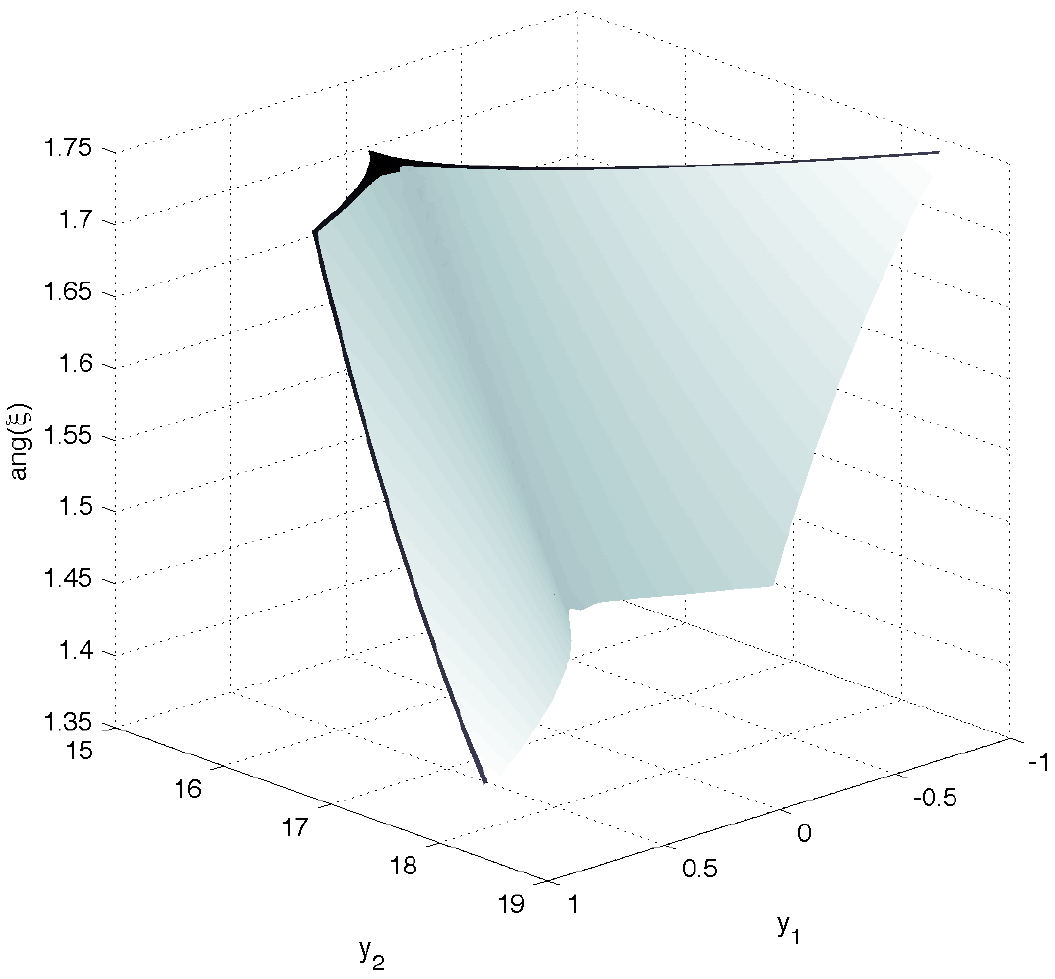}
\caption{\label{fig:caustic} Caustic surfaces  $\Xi(y,\xi)$ (dark gray) and $\check\Xi(y,\xi)$ (light gray) of $\Lambda$ and $\check\Lambda$ corresponding to propagation through a low velocity lens (cf. Section \ref{sec:example}): The singular regions of $\Lambda$ and $\check \Lambda$ do not intersect.}
\end{figure}
%
We index these by $j=1,\ldots,N_i$ and construct a set of diffeomorphisms, $\{ Q_{ij}
\}_{j=1}^{N_i}$, which resolve locally the rank deficiency leading to
coordinates $(y,x_{I_i},\xi_{J_i})$. 
We write
\[
\begin{array}{ccc}
   (y,x_{I_i},\xi_{J_i}) & \stackrel{\kappa_{ij}}{\longrightarrow} &
                 (y,\tilde \xi)
\\
   \uparrow\; \downarrow r & &\;\; \uparrow \;\downarrow \tilde r
\\
   \Lambda \ni r = (y,\eta;x,\xi) \hspace*{1.0cm} & \stackrel{
          C_{Q_{ij}}}{\longrightarrow} & \hspace*{1.0cm}
          (y,\eta;\tilde x,\tilde \xi)
                 = \check r \in \check \Lambda_{ij}
\end{array}
\]
We write $\check O_i$ for the image of $O_i$ under the diffeomorphism
on the level of Lagrangians. Let the matrix $\pdpd{^2 \check
  S_{ij}}{y\partial\tilde \xi}$ in the above be nonsingular on the
open set $\check{U}_{ij}$, and introduce $\check{O}_{ij} =
\check{U}_{ij} \cap \check{O_i} \subset \check \Lambda_{ij}$. This set
corresponds with a set $O_{ij} \subset \Lambda$. We subpartition $O_i =
\cup_{j=1,\ldots,N_i} O_{ij}$. The corresponding partition of unity
now reads
\begin{equation}
   \sum_{i=1}^N \sum_{j=1}^{N_i} \mathit{\Gamma}_{ij}(r) = 1 ,\quad
   \mbox{while}\
   \bar{\mathit{\Gamma}}_{ij}(y,x_{I_i},\xi_{J_i})
               = \mathit{\Gamma}_{ij}(r(y,x_{I_i},\xi_{J_i})) ,\
   j = 1,\ldots,N_i .
\end{equation}
Then $(F \varphi_{\gamma})(y) = \sum_{i=1}^N \sum_{j=1}^{N_i}
(F_{ij} \varphi_{\gamma})(y)$ with
\begin{equation}
   (F_{ij} \varphi_{\gamma})(y)= \int \! \int
           \bar{\mathit{\Gamma}}_{ij}(y,x_{I_i},\xi_{J_i})
   a_i(y,x_{I_i},\xi_{J_i})
   \exp[\ii (S_i(y,x_{I_i},\xi_{J_i})
           - \langle \xi_{J_i},x_{J_i} \rangle]\
   \varphi_{\gamma}(x) 
   \mathrm{d} x \mathrm{d}\xi_{J_i} .
\end{equation}
Inserting the diffeomorphisms, we obtain
\begin{equation}
\label{equ:Fij}
   (F_{ij} \varphi_{\gamma})(y)= \int
   \check{A}_{ij}(y,\tilde{\xi}) \,
   \exp[\ii \check{S}_{ij}(y,\tilde{\xi})]\
   \widehat{Q_{ij}^* \varphi_{\gamma}}(\tilde{\xi})\
   \mathrm{d}\tilde{\xi} .
\end{equation}

\begin{figure}
\centering
\includegraphics[width=\linewidth]{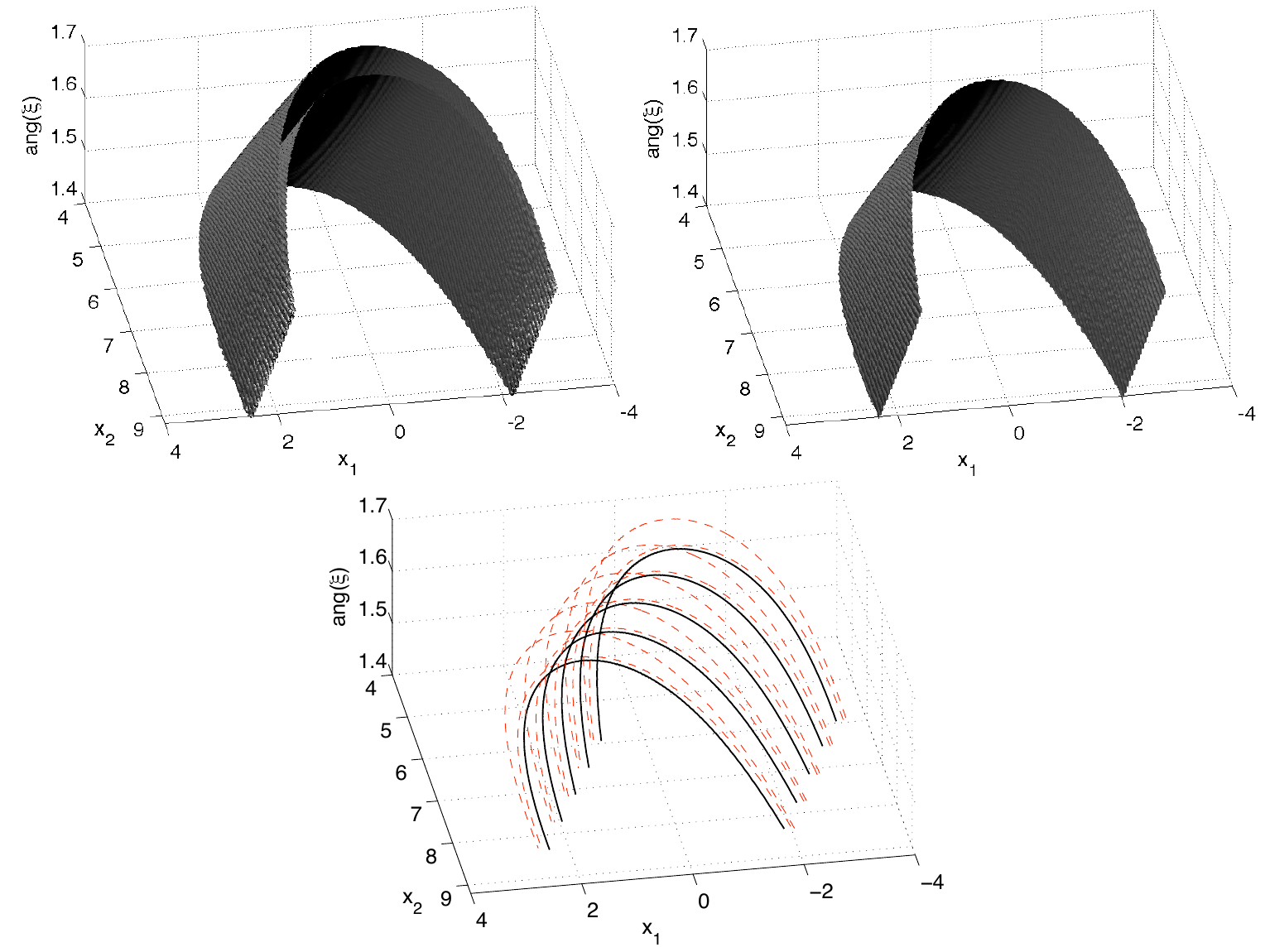}
\caption{\label{fig:iso_win_x_xi} Iso-amplitude surface of the partition functions $\bar\Gamma_i(x,\xi)$, $i=1,3$, associated with $Q_{i}=\mathbb{I}$ (left): the joint admissible set $O_1\cup O_3$ comprises the exterior of the two sheets. 
Iso-amplitude surface of ${\bar\Gamma}_{ij}(x(\tilde x,\tilde\xi),\xi(\tilde x,\tilde\xi))$ for $\xi_0=\pi/2,x_{2,0}=0$ and $\alpha=1$ ($i=2,j=1$) (right): the admissible set $O_{ij}$ contains the region on the back of the sheet. Bottom: boundaries $\partial O_i$, $i=1,3$ (dashed curves) and $\partial O_{ij}$, $i=2,j=1$ (solid curves) of the admissible domains: Clearly, the joint admissible set $O_1\cup O_3 \cup O_{21}$ covers $\Lambda$.}
\end{figure}

The amplitude $\check{A}_{ij}(y,\tilde{\xi})$ and phase function
$\check{S}_{ij}(y,\tilde{\xi}) - \langle \tilde{\xi},\tilde{x}
\rangle$ are obtained by composing $F_{ij}$ with $Q_{ij}^{-1}$ as
Fourier integral operators and changing phase variables. It is possible to treat this composition from a semi-group point of view.
Then, to leading order, we get 
\begin{equation}
\label{equ:amp1}
   \check{A}_{ij}(y,\tilde{\xi}) =
   \bar{\check{\mathit{\Gamma}\ }}_{\!\! ij}(y,\tilde{\xi}) \,
   \check{a}_{ij}(y,\tilde{\xi}) ,
\end{equation}
where
\begin{equation}
   \bar{\check{\mathit{\Gamma}\ }}_{\!\! ij}(y,\tilde{\xi})
       = \check{\mathit{\Gamma} \, }_{\!\! ij}(\check{r}(y,
                       \tilde{\xi})) ,
\end{equation}
in which
\begin{equation}
   \check{\mathit{\Gamma} \, }_{\!\! ij}(\check{r}(r)) =
                                  \mathit{\Gamma}_{ij}(r) .
\end{equation}
Moreover, $\check{a}_{ij}(y,\tilde{\xi})$ can be obtained as
follows. If $\Pi$ is the propagator matrix of the perturbations of
$\chi$, then the propagator matrix of the perturbations of $\check
\chi$ is given by: $\check\Pi_{ij} = \Pi\ \Pi_{Q_{ij}}^{-1}$. Then
\begin{equation}
   \check{a}_{ij}(y,\tilde{\xi}) = \sqrt{1 / \det
   \left(\pdpd{\check S_{ij}(y,\tilde \xi)}{y \partial\tilde \xi}
        \right)^{-1}},
\end{equation}
where $\det
   \left(\pdpd{\check S_{ij}(y,\tilde \xi)}{y \partial\tilde \xi}
        \right)^{-1}$ is obtained as the determinant of the upper-left sub-block of $\check\Pi_{ij}$.
To accommodate a common notation, we set $Q_{ij} = \mathbb{I}$ ($N_i = 1$) if
$I_i = \emptyset$ and write $Q_i$. In the further analysis, we omit the subscripts
${}_{i}$ and ${}_{ij}$ where appropriate.

\begin{figure}
\centering
\includegraphics[width=\linewidth]{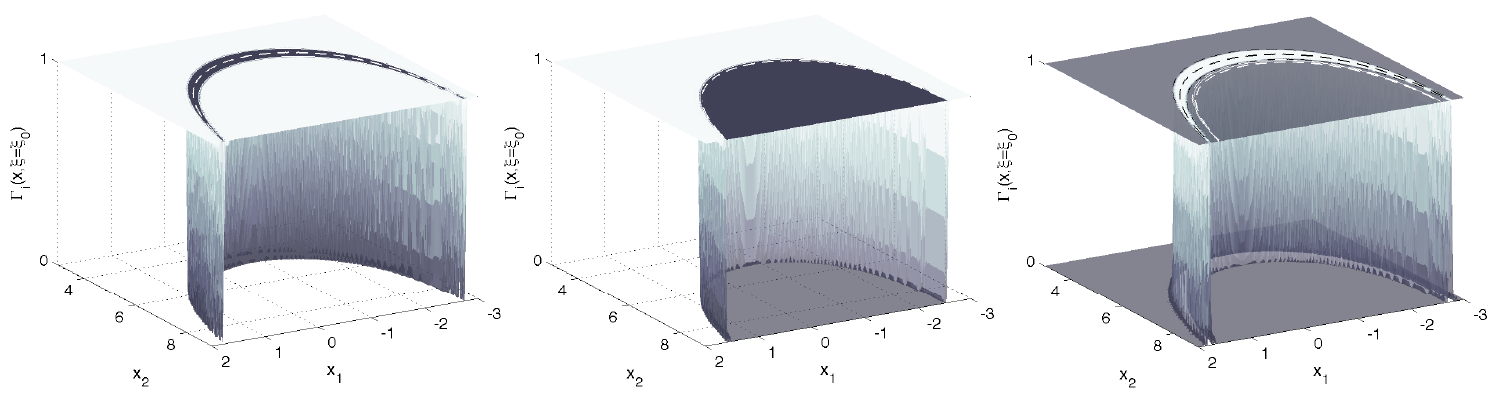}
\caption{\label{fig:Fwin_x_xi} Illustration of joint partition of unity for the partition functions and sets in Fig. \ref{fig:iso_win_x_xi} for $\xi_0=1.67$ fixed: Slice of $\bar\Gamma_i(x,\xi=\xi_0)$ (left), the admissible set $ U_{ij}$
 and the associated partition function ${\bar\Gamma}_{ij}(x(\tilde x,\tilde\xi),\xi(\tilde x,\tilde\xi)=\xi_0)$ (center), and the partition function ${\bar\Gamma}_{ij}(x(\tilde x,\tilde\xi),\xi(\tilde x,\tilde\xi)=1.67)$ for $O_{ij}$ realizing the partition of unity with $\bar\Gamma_i(x,\xi=\xi_0)$.}
\end{figure}

\subsection*{Expansion of the cutoff functions}
To numerically evaluate \eqref{equ:Fij} in reasonable time, we will use separated (in $y$ and $\tilde\xi$) representations of $\check A_{ij}(y,\tilde\xi)$ and $\check S_{ij}(y,\tilde\xi)$ \cite{Andersson2011,Cand`es2007,Cand`es2009}. Such representations can be obtained by restricting the integration over $\tilde\xi$ to domains following a dyadic parabolic decomposition. Here, these will be given by the boxes $B_{\nu,k}(\tilde\xi)$ following the re-decomposition of $Q^*\varphi_{\gamma}$ into wave packets $\varphi_{\tilde\gamma}$, $Q^*\varphi_{\gamma}=\sum_{\tilde\gamma}u_{\tilde\gamma}\varphi_{\tilde\gamma}$. The key novelty is constructing a
separated representation of the partition functions. 

Consider our oscillatory integral in $(y,\tilde{\xi})$ including the
cutoff $\bar{\check{\mathit{\Gamma}\ }} \!\!
(y,\tilde{\xi})$. $\bar{\check{\mathit{\Gamma}\ }} \!\!
(y,\tilde{\xi})$ is homogeneous of degree zero in $\tilde \xi$ and is
a classical smooth symbol (of order $0$). We ``subdivide'' the
integration over $\tilde\xi$. A possible procedure involves obtaining a
(low-rank) separated representation of
$\bar{\check{\mathit{\Gamma}\ }} \!\!  (y,\tilde{\xi})$ on the support
of each relevant box in $\tilde \xi$ \cite{Bao1996,Beylkin2005a,Beylkin2008},
\begin{equation}
\label{equ:expandgamma}
   \bar{\check{\mathit{\Gamma}\ }} \!\!  (y,\tilde{\xi})
   = \sum_{\beta=1}^{J_{\nu,k}}
          \check\Gamma_1^{\beta}(y) \check\Gamma_2^{\beta}(\tilde \xi) ,\quad
   \tilde \xi \in B_{\nu,k} .
\end{equation}
(Basically, this can be obtained using spherical harmonics in view of
the fact that the $\tilde \xi$ is implicitly limited to an annulus.)
One can view this also as windowing the directions of $\tilde \xi$
into subsets (cones) using $\check\Gamma_2^{\beta}(\tilde \xi)$ and then
constructing $\check\Gamma_1^{\beta}(y)$ according to the smallest
admissible set in $y$ for the $\beta$-range of directions.

The oscillatory integral becomes 
    
\begin{equation}
\label{equ:subdiv}
   (F \varphi_{\gamma})(y) = \sum_{\nu,k}
   \sum_{\beta=1}^{J_{\nu,k}}
      \check\Gamma_1^{\beta}(y) \int \check{a}(y,\tilde{\nu})
      \exp[\ii \check{S}(y,\tilde{\xi})] \,
      \check\Gamma_2^{\beta}(\tilde{\xi}) \,
   |\hat{\chi}_{\nu,k}(\tilde{\xi})|^2 \,
          \widehat{Q^* \varphi_{\gamma}}(\tilde \xi)
   \mathrm{d}\tilde{\xi} .
\end{equation}
One can view $\check\Gamma_2^{\beta}(\tilde \xi)
\hat\chi_{\nu,k}(\tilde{\xi})$ as a subdivision of the box
$B_{\nu,k}$. We know that $|J_{\nu,k}| \to 1$ as $k \to \infty$ since
the cone of directions in $B_{\nu,k}$ shrinks as $\sqrt{k}$. Hence,
for large $k$ this does not involve any action.

The procedure allows a subdivision for coarse scales, as long as the
scaling is not affected for large $k$. If the subdivision is too
``coarse'' then parts of the integration will be lost.

\begin{figure}
\centering
\includegraphics[width=0.6\linewidth]{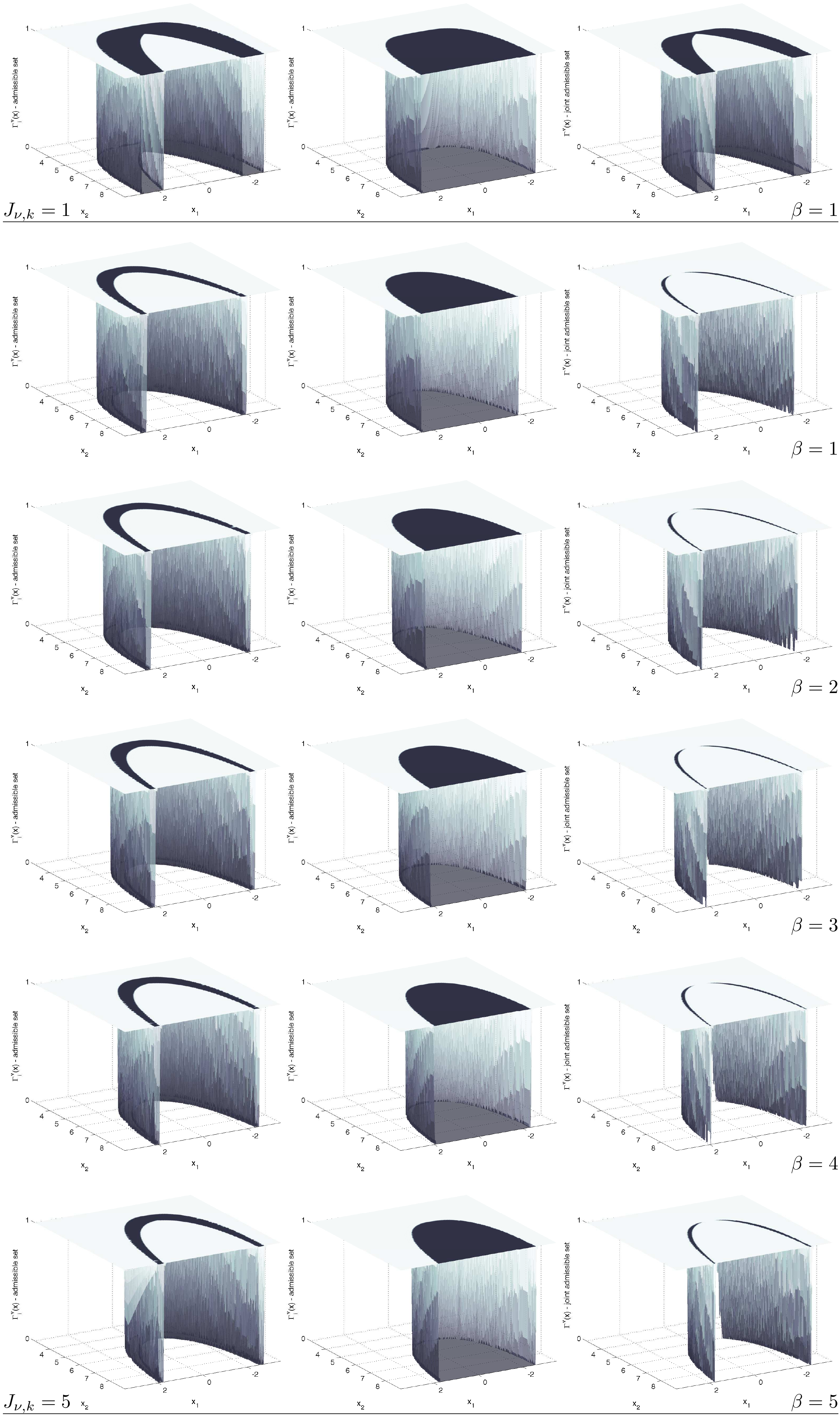}
\caption{\label{fig:F_c_Fwin_x_cc}
Illustration of admissible sets and expansion functions  $\check\Gamma_1^\beta(y(x))$ in \eqref{equ:expandgamma} with  $J_{\nu,k}=1$ (top row) and  $J_{\nu,k}=5$ (bottom rows) for the partition functions  in Fig. \ref{fig:iso_win_x_xi}. Left column: for the partition functions $\bar\Gamma_i(y(x),\xi)$, $i=1,3$. Centre column: for the partition functions ${\bar\Gamma}_{ij}(y(\tilde x,\tilde\xi),\xi(\tilde x,\tilde\xi))$, $i=2,j=1$. Right column: The joint admissible sets induced by the partition of unity for the expansion functions plotted in the left and centre columns.
}
\end{figure}

\begin{figure}
\centering
\includegraphics[width=0.8\linewidth]{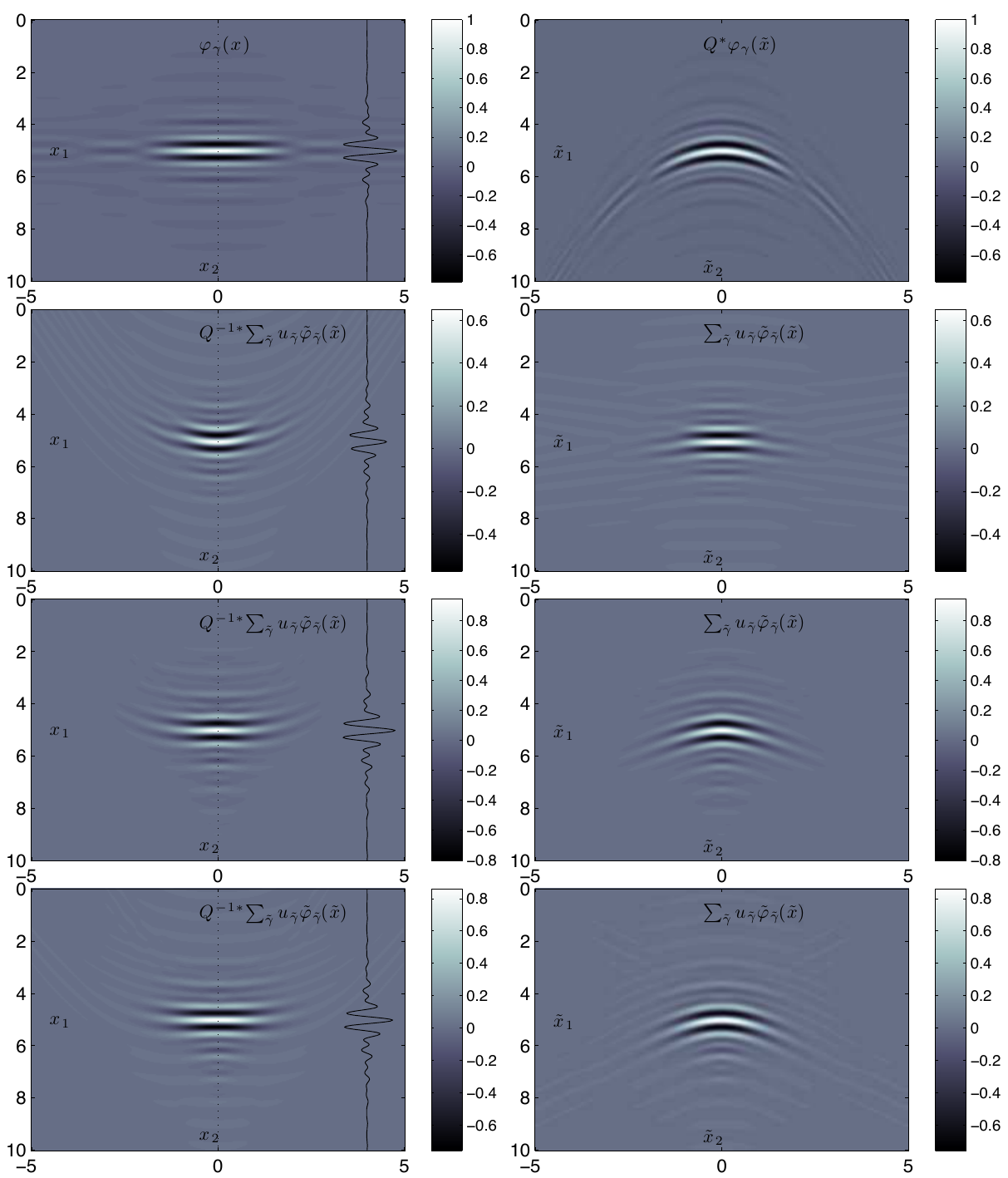}
\caption{\label{fig:recdec}Illustration of diffeomorphism $Q$, $Q^{-1}$ and re-decomposition for a wave packet $\phig(x)$ at frequency scale $k=2$. Top row: $\phig(x)$ (left) and pull-back $Q^*\phig(\tilde x)$ (right). Rows 2 to 4: Re-decompositions $\sum_{\tilde\nu,k}u_{\tilde\gamma}\tilde\varphi_{\tilde\gamma}(\tilde x)$ of $Q^*\phig(\tilde x)$ using $3$, $7$, and $9$ boxes $B_{\tilde\nu,k}$, respectively (right column), and the corresponding image $Q^{*-1}\left(\sum_{\tilde\nu,k}u_{\tilde\gamma}\tilde\varphi_{\tilde\gamma}\right)(x)$ under the action of $Q^{-1}$ (left column). An insufficient number of boxes alters the amplitudes and the minimum phase property of wave packets. Increasing the number of boxes yields satisfactory results in an open neighborhood of $(x_0,\xi_0)$.}
\end{figure}

\clearpage
\section{Computation}
We describe an algorithm for applying Fourier integral operators in the above constructed universal oscillatory integral representation. The global hierarchy of operations is given by the following steps:

\begin{enumerate}
\item {\it Preparation step.\quad} Preparation of universal oscillatory integral representation.
\begin{enumerate}
\item determination of open sets with local coordinates $I_i$, $J_i$ on canonical relation $\Lambda$, inducing $Q_{ij}$
\item construction of cut-off functions ${\check\Gamma}_{ij}$ for the locally singularity resolving diffeomorphisms $Q_{ij}$
\item construction of separated representation for ${\check\Gamma}_{ij}$
\end{enumerate}
\item {\it Evaluation of diffeomorphisms $(Q_{ij}^*\phig)(\tilde x)$.}
\item {\it Evaluation of actions of $(\check F_{ij} (Q_{ij}^*\phig))(y)$.}
\end{enumerate}

Step 3 requires the evaluation  of the action of Fourier integral operators associated with canonical graphs in microlocal standard focal coordinates. The choice of discretization and algorithm for Step 3 induces how computations in Step 1 and 2 are to be organized.
Here, we perform computations in the almost symmetric wave packet transform domain. We make use of the ''box-algorithm'' computation of the action of Fourier integral operators associated with canonical graphs in microlocal standard focal coordinates $(y,\tilde\xi)$ \cite{Andersson2011}. 
The box algorithm is based on the discretization and approximation, to accuracy $\mathcal{O}(2^{-k/2})$, of the action of $\check F_{ij}$ on a wave packet $ \varphi_{j,\tilde \nu,k}(\tilde x)$,
\begin{equation}
\label{equ:boxalgo}
(\check F_{ij} \varphi_{\tilde\gamma})(y)\approx \check A(y,\tilde\nu)\sum_{r=1}^{R} \alpha_{\tilde\nu,k}^{(r)}(y)\sum_{\tilde\xi \in B_{\tilde\nu,k}} e^{\ii \langle T_{\tilde\nu,k}(y),\xi\rangle} |\hat\chi_{\tilde\nu,k}(\tilde\xi)|^2\hat\vartheta_{\tilde\nu,k}^{(r)}(\tilde\xi).
\end{equation}
The procedure relies on truncated Taylor series expansions of $\check S_{ij}(y,\tilde\xi)$ and $\check A(y,\tilde\xi)$ near the microlocal support of $\varphi_{\Tilde\gamma}$, along the $\tilde\nu=\tilde\xi'/|\tilde\xi'|$ axis and in the $\tilde\xi''$ directions perpendicular to the radial $\tilde\nu=\tilde\xi'$ direction. Here, $T_{\tilde\nu,k}(y)$ is the backwards-solution
$$
x(y)=T_{\tilde\nu,k}(y)=\pdpd{\check S_{ij}(y,\tilde\nu)}{\tilde\xi},
$$
and $\alpha_{\tilde\nu,k}^{(r)}(y)$ and $\vartheta_{\tilde\nu,k}^{(r)}(\tilde\xi)$ are functions realizing, on $B_{\tilde\nu,k}$, a separated tensor-product representation of the slowly oscillating kernel appearing in the second-order expansion term of $\check S_{ij}$,
\begin{equation} \label{equ:appro_tensor}
   \exp\left[\ii \frac{1}{2\tilde\xi'}\left\langle\tilde\xi'',\frac{\partial^2\check S_{ij}}{\partial\tilde\xi''^2}(y,\tilde\nu)\;\tilde\xi''\right\rangle \right] B_{\tilde\nu,k}(\tilde\xi)
               \approx \sum_{r=1}^R \alpha_{\tilde\nu,k}^{(r)}(y) \hat\vartheta_{\tilde\nu,k}^{(r)}(\tilde\xi)
\end{equation}
constructed from prolate spheroidal wave functions \cite{Beylkin2005,Slepian1964,Slepian1972,Slepian1978,Xiao2001}.
The number $R$ of expansion terms is controlled by the prescribed accuracy $\varepsilon$ of the tensor product representation.
For a detailed description of the box-algorithm and its implementation, we refer to \cite{Andersson2011}.

Based on this tensor product representation, it is possible to group computations and to evaluate the action of $\check F_{ij}$ in Step 3 for all data wave packets of the same frequency box $B_{\tilde\nu,k}$ at once instead of for each  $\phig$ individually. Consequently, Steps 1 and 2 will also be organized in terms of frequency boxes $B_{\tilde\nu,k}$. We write $u_{\tilde\nu,k}(\tilde x)=\sum_{m}\utgnk\phitgnk(\tilde x)$ for the data portion corresponding to a  frequency box $B_{\tilde\nu,k}$. Starting from data $u(x)$, each of the following steps are repeated for any frequency box $B_{\nu,k}$ of interest.

\subsection{Preparation step}

We begin with determining the sets $O_i$ for the box $B_{\nu,k}$. To this end, we compute
the integral curves $(y(x,\xi),\eta(x,\xi))$ and their perturbations with
respect to initial conditions $(x,\xi)$ and monitor the null space
of the matrix $\pdpd{y}{x}$ following Section
\ref{sec:prop}. For parametrices of evolution equations, this involves solving the system \eqref{equ:ctsystem} and 
evaluating the propagator matrices $\Pi(x,\xi)$ by solving system \eqref{equ:perturbS}. We evaluate the system of differential equations \eqref{equ:perturbS}  in {Fermi} (or ray-centered) coordinates, in which the potential rank deficiencies of the upper left subblock $W_1(x,\xi)$ appear explicitly as zero entries in the corresponding row(s) and column(s) \cite{Cerveny2001}. The submanifolds $\Sigma_{(x,\xi)}$ on which $W_1$ is singular separate and define the sets $O_i$. These computations on $\Lambda_F$ are performed by discretizing the set of orientations $\nu=\xi/|\xi|$ covering the frequency box $B_{\nu,k}$ with resolution $\delta_\nu$ and the set in $x$ for which $u_{\nu,k}(x)$ has non-zero energy with resolution $\delta_x$. 

Then, for each set $O_i$, we detect $\check U_{ij}$ (and consequently $\check O_{ij}$) in
a similar way, as the set on which the upper left sub-block $\check W_{1,ij}$ of
$\check\Pi_{ij}=\Pi\;\Pi^{-1}_{Q_{ij}}$ has full rank. Here
$\Pi^{-1}_{Q_{ij}}(\tilde x,\tilde\xi)$ is given by
\eqref{equ:propQ1}.
The operators $Q_{ij}$ are chosen such that $\{\check O_{ij}\}_{j=1}^{N_i}$ overlappingly cover the set of singularities. For fixed $\alpha_i$, this induces a discrete set $\{x_0^j\}_{j=1}^{N_i}$.

We then proceed with the construction of the partition of unity. Since the partition functions enter the computation as pseudodifferential cutoffs in the construction of the amplitude (cf. \eqref{equ:amp1}), requiring the backwards solutions $\tilde x(y,\tilde \xi)$ (compare (\ref{equ:d2S1}--\ref{equ:amp0})), we perform our numerical construction in coordinates $(\tilde x,\tilde\xi)$. We obtain $\bar{\check\Gamma}_{ij}(y,\tilde\xi)$ upon substituting $y=y(\tilde x,\tilde\xi)$ implied by the canonical relation $\check \chi_{ij}$. For the construction of the partition functions $\check\Gamma_{ij}$, we choose double-exponential cutoffs of the form
$$
\exp(-\exp(d(\tilde x,\tilde\xi))
$$
mimicking a $C_0^\infty$ cutoff, with appropriate normalization and truncated to precision $\varepsilon$.  Here $d(\tilde x,\tilde\xi)$ is a function measuring the distance of the point $(\tilde x,\tilde\xi)$ from the boundary $\partial \bar{\check U}_{ij}$ of the set $\bar{\check U}_{ij}(\tilde x,\tilde\xi)$. The partition of unity is then formed by weighting $\bar{\check \Gamma}_{ij}(\tilde x,\tilde\xi)$ on the overlaps of the sets $\bar{\check U}_{ij}(\tilde x,\tilde\xi)$ such that $\sum_{ij}{\check \Gamma}_{ij}(\check r(\tilde x,\tilde\xi))=1$.

Finally, we construct the separated representations of $\check\Gamma_{ij}$ (cf. \eqref{equ:subdiv}) in $(\tilde x,\tilde\xi)$ coordinates by windowing the directions of $\tilde\xi$ into subsets using $\check\Gamma_2^\beta(\tilde\xi)$, realizing a subdivision into $\tilde \xi$ cones. This subdivision is performed for each frequency box $B_{\tilde\nu,k}$. 

\subsection{Evaluation of  diffeomorphisms}

We evaluate each of the operators $Q_{ij}$ in the Fourier domain. 
This choice is guided by the property  $\sum_m\ugnk\hat{\varphi}_{m,\nu,k}(\xi) = \hat{u}(\xi)
\hat{\beta}_{\nu,k}(\xi) \hat{\chi}_{\nu,k}(\xi)$ of the discrete almost symmetric wave packet transform \cite{Duchkov2009a}, which enables the fast evaluation of the Fourier transform of the data at a set of frequency points $\xi_l^{\nu,k}$ limited to the box $B_{\nu,k}$. We obtain $(Q_{ij}^*\phig)(\tilde x)$ at once for all $\phig(x)$ belonging to the frequency box $B_{\nu,k}$  by evaluation of their adjoint unequally spaced FFT \cite{Dutt1993,Dutt1995}, $\mathcal{F}^{US\,*}_{ \xi\;\rightarrow\;x}$, at points $x(\tilde x)=Q_{ij}^{-1}(\tilde x)$,
$$
\check u_{\nu,k}^{ij}(\tilde x)=\sum_{m}\ugnk(Q_{ij}^*\phignk)(\tilde x)=\mathcal{F}^{US\,*}_{ \xi=\xi_l^{\nu,k}\;\rightarrow\;x(\tilde x)}\big[\hat{u}(\xi)
\hat{\beta}_{\nu,k}(\xi) \hat{\chi}_{\nu,k}(\xi)\big]
$$
In preparation for the evaluation of $(\check F_{ij}\check u_{\nu,k})(y)$, we compute the discrete almost symmetric wave packet transform of the pullback $\check u^{ij}_{\nu,k}(\tilde x)$, yielding its wave packet coefficients $u_{\mathbb{j},\tilde\nu,k}^{ij}$. 

\subsection{Evaluation of the actions of $F_{ij}$}

At this stage, we are ready to evaluate the action 
$(F_{ij}u_{\nu,k})(y) = \sum_m\ugnk(F_{ij}\phignk)(y)$ (cp. \eqref{equ:Fij}) by evaluation of $(\check F_{ij}\check u^{ij}_{\nu,k})(y)$
using the box algorithm (cf. \eqref{equ:boxalgo}).
Note that numerically significant coefficients $u^{ij}_{\mathbb{j},\tilde\nu,k}$ of the pull-back $\check u^{ij}_{\nu,k}(\tilde x)$ are contained in a small set of boxes $B_{\tilde\nu,k}$ neighboring the direction $\nu=\xi_0/|\xi_0|$. 
We further subdivide each of these boxes according to the separated representation of $\check\Gamma_{ij}$.
Then, we apply the box algorithm to each subdivision, indexed by triples $(\beta,\tilde\nu,k)$,  $\beta=1,\dots,J_{\tilde\nu,k}$.
Here, the Taylor series expansion of the generating function $\check S_{ij}(y,\tilde\xi)$ is constructed about the central $\tilde\xi$ direction within the support of $\Gamma_2^{\beta}(\tilde \xi) \hat\chi_{\tilde\nu,k}(\tilde{\xi})$, accounting for the induced subdivision of the box $B_{\tilde\nu,k}$. 
Note that sub-dividing into $\tilde\xi$ cones results in a reduction of the range of $\tilde\xi$ orientations in each element $(\beta,\tilde\nu,k)$ of the subdivision, as compared to the $\tilde \xi$ range contained in $B_{\tilde\nu,k}$. This reduces the number $R$ of expansion terms in \eqref{equ:boxalgo} and effectively counter-balances the increase  by a factor $J_{\tilde\nu,k}$, evoked by the separated representation of $\check\Gamma_{ij}$, of the number of times  the box-algorithm has to be applied.


\subsection*{Operator hierarchy}

The operators $F_{ij}$ for which $Q_{ij}=\mathbb{I}$, $F_{ij}^\mathbb{(I)}$ say, are directly associated with the canonical relation $\Lambda_F$ and involve only computations on $\Lambda_F$.
In the algorithm, we reflect this physical hierarchy of the operators $F_{ij}$ in the construction of the partition of unity. First, we construct a partition of unity for these hierarchically higher operators. Then, we construct a joint partition of the remaining operators on the sets which are not covered by the sets for which $Q_{ij}=\mathbb{I}$.

\subsection*{Re-decomposition}

Starting from a single box $B_{\nu,k}$ and applying $Q_{ij}$ to $u_{\nu,k}(x)$, re-decomposition of $\check u^{ij}_{\nu,k}(\tilde x)$ results in a set of boxes $B_{\tilde\nu,k}$ yielding numerically non-zero contribution to the solution.
The number of boxes entering the computation is directly proportional to the computational cost of the algorithm. In applications, we therefore aim at keeping this number small and consider only a subset of boxes, yielding the most significant contributions. We choose this subset such that on an open neighborhood of $(x_0,\xi_0)$
$$Q_{ij}^{-1}Q_{ij}\approx \mathbb{I}$$
to precision  $\epsilon$.
We can estimate the energy loss induced by the restriction to subsets of $B_{\tilde\nu,k}$ and re-normalize the solution. We illustrate the impact of choices of subsets containing different numbers of boxes on the numerical accuracy of the diffeomorphic identity in Fig. \ref{fig:recdec}.


Furthermore, the re-decomposition of $\check u^{ij}_{\nu,k}(\tilde x)$ yields in general, under the action of $Q_{ij}^{-1}$, $\xi$-values outside the set $B_{\nu,k}$, $\xi(x,\tilde\xi)\supset B_{\nu,k}$. We monitor $\xi(x,\tilde\xi)$ and do not consider their contribution in our computation if $|\hat\chi_{\nu,k}(\xi(x,\tilde\xi))|$ is below the threshold $\varepsilon$.

The free parameters of the procedure are summarized in Table \ref{tab:param}.

\begin{table}
\centering
\begin{tabular}{ll}
$\delta_\nu$, $\delta_x$ & discretization steps in computations on $\Lambda_F$ and $\check \Lambda_{ij}$\\
$\alpha_i$ & free parameter of operators $Q_{ij}$, inducing the discrete sets of diffeomorphisms $\{Q_{ij}\}_{j=1}^{N_i}$\\
$J_{\nu,k}$ & number of expansion terms in separated representation of cutoff functions $\check \Gamma(y,\tilde\xi)$\\
 $\epsilon$ & precision of approximate re-decomposition  of $\check u^{ij}_{\nu,k}(\tilde x)$\\
 $\varepsilon$ & accuracy of the tensor-product representation in the box algorithm
\end{tabular}
\caption{\label{tab:param}Table of user defined parameters.}
\end{table}

\begin{figure}
\centering
\includegraphics[width=0.66\linewidth]{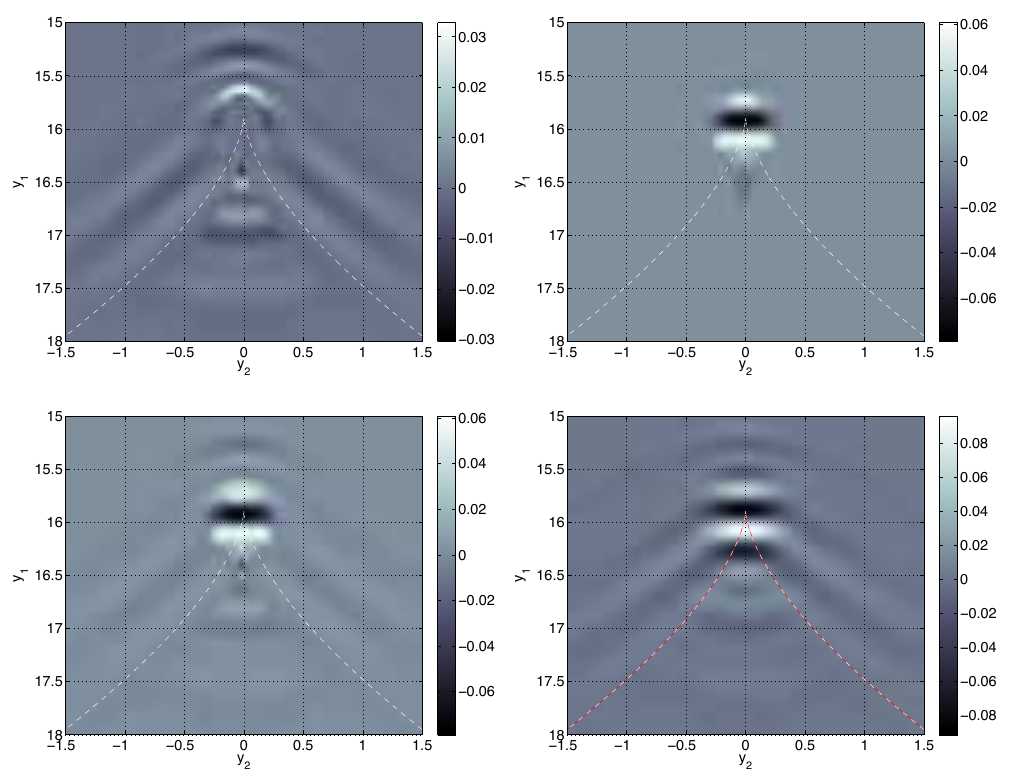}
\caption{\label{fig:3WP0}
Illustration of operator action on a wave packet $\phig(x)$ at frequency scale $k=2$: Contribution of operators $F_i$ ($i=1,3$) associated with $Q_{i}=\mathbb{I}$ (top left), contribution of operator $F_{ij}$ ($i=2,j=1$) with diffeomorphism parameters $(\xi_0=\pi/2,x_{2,0}=0,\alpha=1)$, resolving the singularity in the tip of the caustic (top right), and joint action of $F_i$ and $F_{ij}$ (bottom left). Time domain finite difference reference (bottom right). In the operator computation, we consider 9 boxes $B_{\tilde\nu,k}$ and a separated representation with $J_{\nu,k}=1$ term.
}
\end{figure}
\begin{figure}
\centering
\includegraphics[width=\linewidth]{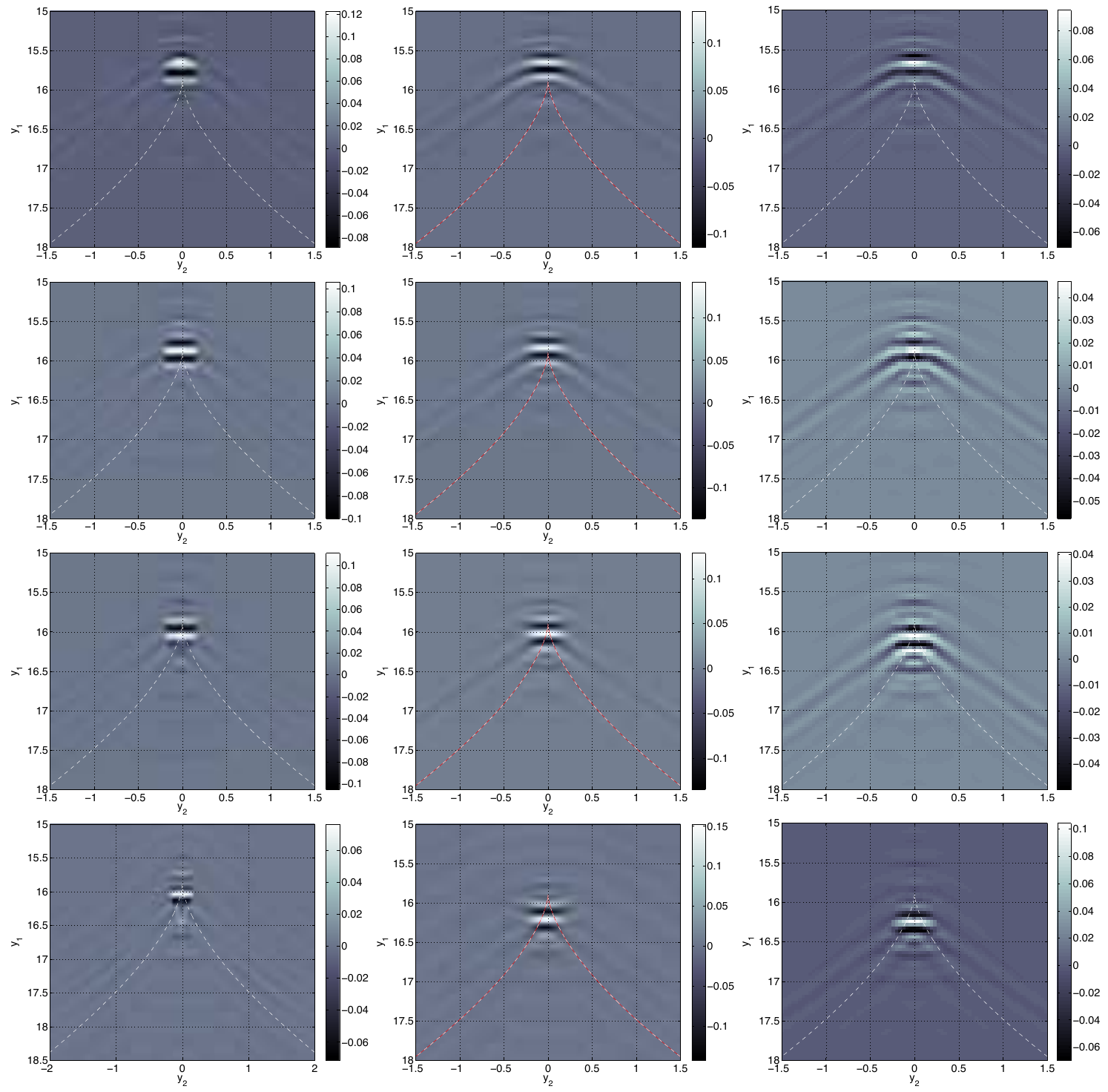}
\caption{\label{fig:4WP4}
Left column: Joint contribution of the operators $F_i$ and $F_{ij}$ acting on a wave packets $\phig(x)$ at frequency scale $k=3$ (compare Fig. \ref{fig:3WP0} (bottom left) for a wave packet at frequency scale $k=2$). Center column: Time domain finite difference reference. Right column: The equivalent to the left column when using a separated representation with $J_{\nu,k}=11$ terms (note that for computational reasons, only $1$ box $B_{\tilde\nu,k}$ has been used in the numerical evaluation of $F_i$ and $F_{ij}$ with $J_{\nu,k}=11$ terms).
}
\end{figure}

\section{Numerical example}
\label{sec:example}

We numerically illustrate our algorithm for the evaluation of the
action of Fourier integral operators associated with evolution
equations. We consider wave evolution under the half-wave equation,
that is, the initial value problem \eqref{equ:evo:equ} with symbol
$$
P(x,\xi)=\sqrt{c(x)^2||\xi||^2},
$$
in $n=2$ dimensions. Here $c(x)$ stands for the medium velocity. 

\subsection*{Heterogeneous, isotropic model}
We choose a heterogeneous velocity model
$$
c(x)=c_0+\kappa \exp(-|x-x_0|^2/\sigma^2),
$$
containing a low velocity lens, with parameters $c_0=2km/s$, $\kappa=-0.4km/s$, $\sigma=3km$, and $x_0=(0,14)km$.
As the initial data, we choose horizontal wave packets at frequency scale $k=2$ and $k=3$, respectively, in the vicinity of the point $x'=(0,5)km$. We set $t_0=0$ and fix the evolution time to $T=7s$. With this choice of parameters, most of the energy of the solution is concentrated near a cusp-type caustic. We illustrate the induced sets $O_i$ and the joint partition of unity $\Gamma_i$ in Fig. \ref{fig:iso_win_x_xi} and \ref{fig:Fwin_x_xi}.

\subsection*{Operator factorization}
We partition the Lagrangian $\Lambda$ into three sets $O_i$, $i=\{1,2,3\}$. The sets $i=\{1,3\}$ are separated by the caustic. For these sets, we can choose coordinates $(y,\xi)$, hence $Q_i=\mathbb{I}$. The set $i=2$ contains the caustic.
For illustration purposes, in the factorization $F_{ij}$ of $F_i$ for $i=2$, we choose to compute the operator $j=1$, which   resolves the singularity in an open neighborhood of the point indicated by a black dot on the Lagrangian plotted in Fig. \ref{fig:lagrangian}. This neighborhood contains the cusp of the caustic. Furthermore, we limit our separated representation to one term, $J_{\nu,k}=1$ (for the corresponding admissible sets and partition functions, see Fig. \ref{fig:F_c_Fwin_x_cc} (left column)). 
We restrict the computation of $F_{ij}$ for the initial data at frequency scale $k=2$ ($k=3$) to $9$ ($11$) boxes $B_{\tilde\nu,k}$ neighboring the $\nu$ direction, respectively.

\subsection*{Results}
In Fig. \ref{fig:3WP0}, we plot the contributions of the different components in the factorization of the propagator acting on a single horizontal wave packet at frequency scale $k=2$, and compare to a time domain finite difference computation. The support of the wave packet within the joint admissible set of the chosen factorization is mostly covered by the set $O_{ij}$,  such that most of its energy is contributed by the operator $F_{ij}$, for which $Q_{ij}\neq\mathbb{I}$. 

We observe that in the joint admissible set, our algorithm has effectively removed the singularity. We note that the phase of the operator computation matches the phase of the finite difference reference. This includes the KMAH index, which is best observed for operator $F_3$, which exclusively contributes to the region beyond the caustic (cf. Fig. \ref{fig:3WP0}, top left). Furthermore, note that the amplitude obtained by our algorithm is slightly weaker than the true amplitude. This is consistent with the observations and discussion following Fig. \ref{fig:recdec} and results from the energy leakage induced by restricting the number of boxes in the re-decomposition step following the application of $Q$. We can compensate and re-normalize the amplitude by monitoring the energy loss resulting from the restriction (in Fig. \ref{fig:3WP0}, we have not re-normalized the amplitudes).
Finally, we note that our algorithm yields the correct result in an open neighborhood in the vicinity of the tip of the caustic, for which we have designed the operator $F_{ij}$. In consistency with this fact, it is ineffective for yielding the image of the entire wave packet which, at this low frequency scale, has support extending beyond the admissible set of the operator factors we compute. 

\medskip

These observations are further illustrated in Fig. \ref{fig:4WP4} (left column), where we plot the contributions of the different components in the factorization of the propagator acting on horizontal wave packets, at higher frequency scale $k=3$, centered at $4$ different locations in the vicinity of the caustic tip. Results of a time domain finite difference reference computation are plotted in Fig. \ref{fig:4WP4} (center column). With these initial data, we explore the open neighborhood about the point for which the operator composition with $Q_{ij}$ resolves the singularity. Indeed, at this frequency scale, we can obtain the image of an entire wave packet with only one operator factor $F_{ij}$ (cf. Fig. \ref{fig:4WP4} (second row)). For the wave packet located slightly further above the tip of the caustic (top row), we observe a phase artifact in the region of overlap of $O_{i=1}$ and $O_{ij}$, which can be explained as follows: The restriction of the separated representation for $F_i$ to one term only induces that the computation of the geometry (bi-characteristics) for the entire box $B_{\nu,k}$ is exclusively based on one single direction $\nu$. This results in inaccuracies in regions close to the caustics where slight perturbations in $\xi$ yield large variations in $y$. Furthermore, as discussed above, wave packets exploring the regions beyond the tip of the caustics eventually start to leave the admissible set for $F_{ij}$ (third and bottom line).

We note that both for removing the phase artifact of $F_i$ close to
the caustic, and for enlarging the admissible set, it is necessary to
increase the number of terms $J_{\nu,k}$ in the separated
representation \eqref{equ:expandgamma} (compare
Fig. \ref{fig:F_c_Fwin_x_cc}). This is illustrated in
Fig. \ref{fig:4WP4} (right column), where the joint contributions of
operators $F_i$ and $F_{ij}$ with $J_{\nu,k}=11$ terms in the
separated representation are plotted. Here, the expansion functions
$\check\Gamma_2^\beta(\tilde\xi)$ are constructed as as cones in
$\tilde\xi$ with a squared cosine cutoff window. For practical reasons
and illustration purpose, only one single frequency box
$B_{\tilde\nu,k}$ has been used in the computation
(cf. Fig. \ref{fig:recdec}). While this restriction to only one
frequency box affects the amplitudes and the phases in the tails of
the wave packet, the separated representation remains nonetheless
effective in resolving the issues observed above: the admissible set
is extended beyond the caustic and the inaccuracies in the regions of
overlap of sets $O_i$ and $O_{ij}$ as well as in the regions close to
the caustics are considerably reduced.

\section{Discussion}

We developed an algorithm for the evaluation of the action of Fourier
integral operators through their factorization into operators with a
universal oscillatory integral representation, enabled by the
construction of appropriately chosen diffeomorphisms. The algorithm
comprises a preparatory geometrical step in which open sets are
detected on the canonical relation for which specific focal
coordinates are admissible. This covering with open sets induces a
pseudodifferential partition of unity. Then, for each term of this
partition, we apply a factorization of the associated operators using
diffeomorphisms reflecting the rank deficiency and resolving the
singularity in the set. This factorization admits a parametrization of
the canonical graph in universal $(y,\tilde\xi)$ coordinate pairs and
enables the application of our previously developed box algorithm,
following the dyadic parabolic decomposition of phase space, for
numerical computations. Hence, our algorithm enables the discrete wave
packet based computation of the action of Fourier integral operators
globally, including in the vicinity of caustics. This wave packet
description is valid on the entire canonical relation. It can now
enter procedures aiming at the iterative refinement of approximate
solutions, and drive the construction of weak solutions via Volterra
kernels \cite{AdHSU,Hoop2010}.
 
An alternative approach for obtaining solutions in the vicinity of
caustics has been proposed previously
\cite{Andersson2011,Kumano-go1978,Rousseau2006} for the special case of Fourier integral operators corresponding to
parametrices of evolution equations for isotropic media. 
It consist in a
re-decomposition strategy following a multi-product representation of
the propagator. Here, we avoid the re-decompositions and operator
compositions following the discretization of the evolution parameter,
reminiscent of a stepping procedure. What is more, our construction is
not restricted to parametrices of evolution equations, but is valid
for the general class of Fourier integral operators associated with
canonical graphs, allowing for anisotropy. The cost of the algorithm
resides in the construction and application of the separated
representation of the pseudodifferential partition of unity.

\end{document}